\documentclass[11pt]{article}

\usepackage{placeins}
\usepackage{graphicx,verbatim,color}
\usepackage{enumitem}
\usepackage{natbib}
\usepackage{amsmath}
\usepackage{amssymb}
 \usepackage{multirow}
 \usepackage{mathrsfs}
\usepackage{amsfonts}
\usepackage{algorithm}
\usepackage{algorithmic}
\usepackage{graphicx, subfig}
\usepackage{epsfig}
\usepackage{epstopdf}

\usepackage{amsthm}\usepackage{dsfont}\usepackage{array}\usepackage{mathrsfs}\usepackage{comment}
\usepackage{hyperref}
\usepackage{breakurl}

\usepackage[utf8]{inputenc} 
\usepackage[T1]{fontenc}    
\usepackage{hyperref}       
\usepackage{url}            
\usepackage{booktabs}       
\usepackage{amsfonts}       
\usepackage{nicefrac}       
\usepackage{microtype}  

\usepackage{bm}

\newcommand{\beq}{\begin{equation}}
\newcommand{\eeq}{\end{equation}}

\newcommand{\lb}{\lambda}

\definecolor{darkred}{RGB}{150,0,0}

\definecolor{emmanuel}{RGB}{255,127,0}

\newcommand{\argmin}{\mathop{\mathrm{argmin}}}
\newcommand{\argmax}{\mathop{\mathrm{argmax}}}

\newcommand{\R}{\mathbb{R}}

\newcommand{\ones}{\mathbf{1}}
\newcommand{\zeros}{\mathbf{0}}
\newcommand{\ra}{\rightarrow}
\newcommand{\la}{\leftarrow}
\newcommand{\Ra}{\Rightarrow}

\title{A short note on solving box inequality and linear equality constrained optimization problem}
\author{Yue Sun}
\begin{document}

\theoremstyle{plain}\newtheorem{lemma}{\textbf{Lemma}}\newtheorem{condition}{\textbf{Condition}}\newtheorem{theorem}{\textbf{Theorem}}\newtheorem{corollary}{\textbf{Corollary}}\newtheorem{example}{\textbf{Example}}\newtheorem{definition}{\textbf{Definition}}\newtheorem{conjecture}{\textbf{Conjecture}}
\newtheorem{claim}{\textbf{Claim}}\newtheorem{remark}{\textbf{Remark}}\newtheorem{question}{\textbf{Question}}
\theoremstyle{definition}


\maketitle
\begin{abstract}
    This writeup discusses a special kind of convex constrained optimization problem, whose constraints consist of box inequalities and linear equalities. For this problem, in addition to general optimization algorithms such as exact penalty algorithm and interior point algorithm (\cite{boyd2004convex,BertsekasDimitriP.1999Np}), there is a simple iterative algorithm that is simple to implement, which is favored by machine learning practitioners. 
\end{abstract}
\section{Introduction}
It is known for long that iterative algorithms that follows the gradient is guaranteed to converge to the global optimum of a convex optimization problem. The algorithms for both unconstrained and constrained problems are well studied. For unconstrained problem, one can simply proceed by descent algorithm based on gradient, whereas barrier methods are required for constrained problems such as exact penalty algorithm and interior point algorithm (\cite{boyd2004convex,BertsekasDimitriP.1999Np}). 

Though they are both polynomial time algorithms, the latter is not preferred by practitioners since they are not easy to code. It is hoped that projected gradient algorithm (\cite{bubeck2015convex}) can be applied, but projection is not computable for all kinds of convex sets. This writeup argues that, the projection onto a special kind of constraints, which consist of only box inequalities and linear equalities, is easy to do.

The motivation of such constraints comes from two problems. The first is probability simplex, which is defined as 
\begin{align}\label{eq:prob_simp}
    \{x\in\R^n \ |\ x\ge 0,\ \ones^Tx=1\}.
\end{align}
Thus here $x$ is an $n$-point probability distribution. This originates from the authors work \cite{sun2020online}, and is more widely used in policy training in reinforcement learning problems. Another motivation is \cite{chen2020learning}, which indicates that, frustrated by the obscurity of constrained optimization algorithms, practitioners are forced to switch to a neural network to regress the solution of the problem, though their complexity are both polynomial. Apparently, neural networks may not converge to the optimum, and generally it does not work. Though \cite{chen2020learning} proposes a novel duality trick aiming for the box inequalities and linear equalities structure, unfortunately the recovery is not guaranteed for sure either. Thus, this writeup urges a simple-to-implement algorithm for this kind of problem, which both converges theoretically with polynomial rate as general optimization algorithms, and is easy to implement. 

\section{Simple case: optimization problem constrained on probability simplex}
Let $f(x):\R^n\rightarrow \R$ be a smooth convex function. We consider minimizing $f$ on probability simplex (\ref{eq:prob_simp}), which writes
\begin{align}\label{eq:opt_prob_simp}
    \min\limits_{x\ge0, \ones^Tx=1} \ f(x).
\end{align}
Let $\eta$ be the step size projected gradient descent algorithm writes the iteration as
\begin{align}
    x_{t+1} = \argmin\limits_{x\ge0, \ones^Tx=1}\ \|x - (x_t - \eta \nabla f(x_t))\|_2.
\end{align}
For simplicity of notation, we denote
$y = x_t - \eta \nabla f(x_t)$, and consider
\begin{align}
    \argmin\limits_{x\ge0, \ones^Tx=1}\ \frac{1}{2}\|x - y\|^2.
\end{align}
We write the Lagrangian
\begin{equation}
    L(x, \lb, \mu) = \frac{1}{2} \|x-y\|^2 - \lb^Tx + \mu (\ones^Tx - 1).\vspace{-.4em}
\end{equation}
The KKT condition (fixing $\mu$)  is
\begin{align}
    \nabla_x L(x,\lb) = x-y-\lb+\mu\ones=0;\
    \lb\ge0; \ x\ge0;\ \lb_ix_i=0, \forall i.
\end{align}
The optimal $x$ is $x_i^* = \max(y_i - \mu^*, 0)$.

So we solve 
\begin{align}
    \max_\mu\ \frac{1}{2} \|\max(y-\mu\ones, \zeros)-y\|^2 + \mu (\ones^T\max(y-\mu\ones, \zeros) - 1) \label{eq:dual}
\end{align}
to get $\mu^*$, and $x^* = \max(y-\mu^*\ones, \zeros)$. Since (\ref{eq:dual}) is 1-dimensional convex optimization problem, we can solve it by binary search, which converges in logarithmic time.

\begin{algorithm}[htbp!]
\caption{Algorithm for solving (\ref{eq:opt_prob_simp})}
\begin{algorithmic}
\REQUIRE{Smooth function $f$, step size $\eta$, initial point $x_0$.}\\
\STATE{$t\la 0$.}
\WHILE{not converge}
\STATE{$y\la x_t - \eta \nabla f(x_t)$.}
\STATE{$\mu \la \argmax_\mu\ \frac{1}{2} \|\max(y-\mu\ones, \zeros)-y\|^2 + \mu (\ones^T\max(y-\mu\ones, \zeros) - 1)$.}
\STATE{$x_{t+1} = \max(y-\mu\ones, \zeros)$.}
\STATE{$t\la t+1$.}
\ENDWHILE
\RETURN{$x_t$.}
\end{algorithmic}
\end{algorithm}

\section{Application: network with independent edge flows}
Consider the \emph{network with independent edge flows} problem in \cite{chen2020learning}. It can be simplified as a linear programming:
\begin{align}\label{eq:opt_net}
    \min\limits_{u\le x\le v, Ax=b} \ c^T x.
\end{align}
where $u,v,A,b,c$ are fixed constants and $x\in\R^n$ is the optimization variable. It is more difficult since the equality constraint is not $1$ dimensional. If its dimension is higher than constant, then binary search type algorithms will end up exponential time. However, since the objective is linear, we can use a similar trick with small modification.

First, consider
\begin{align}\label{eq:net}
    \min\limits_{u\le x\le v, Ax=b} \ \frac{1}{2}\|x-y\|^2.
\end{align}
We write the Lagrangian
\begin{equation}
    L(x, \lb, \mu) = \frac{1}{2} \|x-y\|^2 + \lb_1^T(u-x) + \lb_2^T(x-v) + \mu^T (Ax-b).
\end{equation}
The KKT condition is
\begin{gather}
    \nabla_x L(x,\lb_1,\lb_2) = x-y-\lb_1+\lb_2+A^T\mu=0;\\
    \lb_1\ge0; \ \lb_2\ge 0;\\
    u\le x\le v;\\
    \lb_{1,i}(x-u)_i = 0;\ \lb_{2,i}(x-v)_i=0,\ \forall i.
\end{gather}
Denote the optimal $\mu$ as $\mu^*$, we have \begin{align}
    (x,\lb_1,\lb_2)_i^* &=
    \begin{cases}
    ((y-A^T\mu^*)_i,0,0),&\ u_i\le y_i - A^T\mu^*\le v_i;  \\
    (u_i,u_i-(y - A^T\mu^*)_i,0),&\ y_i - A^T\mu^*\le u_i;\\
    (v_i,0, (y - A^T\mu^*)_i-v_i),&\ y_i - A^T\mu^*\ge v_i.
    \end{cases}\\
    &= \phi(\mu^*).\label{eq:phimu}
\end{align}
So the optimal $(x,\lb_1,\lb_2)$ tuple is a function of $\mu$, denoted by $\phi$. Thus we can solve
\begin{align}\label{eq:mu_prob}
    \max_\mu\ L(x,\lb_1,\lb_2, \mu) = \max_\mu\ L(\phi(\mu))
\end{align}
to get $\mu^*$, and return corresponding $x^*$ via $\phi$ function. Note that (\ref{eq:mu_prob}) is not smooth, but zeroth order methods work widely well in unconstrained problem solving by machine learning trainers.

However, (\ref{eq:net}) is not (\ref{eq:opt_net}), and if we use an outer loop for (\ref{eq:opt_net}) and an inner loop for (\ref{eq:net}), the one more layer costs higher computational complexity (though it's acceptable in a lot of regimes).

Fortunately, the linear objective can simplify it in the following way. We state the following theorem:
\begin{theorem}\label{thm:bound_opt}
Denote $S:= \{x\in\R^n\ |\ u\le x\le v,\ Ax=b\}$. The optimizer of (\ref{eq:opt_net}) coincides with
\begin{align}
    \hat x = \mathop{\mathrm{lim}}\limits_{t\ra+\infty}\argmin\limits_{u\le x\le v, Ax=b} \ \frac{1}{2}\|x+tc\|^2.
\end{align}
Moreover, let the radius of $S$ be $r$ and $\max_{x\in S} \|x\| = R$, denote
\begin{align}
    \hat x = \argmin\limits_{u\le x\le v, Ax=b} \ \frac{1}{2}\|x+tc\|^2.
\end{align}
then 
\begin{align*}
    c^T (\hat x - x) \le 4Rr/t.
\end{align*}
\end{theorem}
\textbf{Proof.} We prove by contradiction, say, suppose
\begin{align}
   c^T\hat x \ne c^T x^*.   
\end{align}
WLOG let $c\ne 0$. Let the true optimizer of (\ref{eq:opt_net}) be $x^*$, from its optimality, we have that
\begin{align}
    c^T(x - x^*)\ge 0,\ \forall x\in S.
\end{align}
Thus $\delta:=c^T(\hat x - x^*)>0$ because $\hat x\in S$. From optimality of $\hat x$ we have
\begin{align}
    &\quad \|\hat x+tc\|^2 \le \|x^*+tc\|^2\\
    &\Ra \|\hat x\|^2 + 2tc^T\hat x + t^2\|c\|^2 \le  \|x^*\|^2 + 2tc^Tx^* + t^2\|c\|^2\\
    &\Ra 2tc^T(\hat x - x^*) \le -(\hat x+x^*)^T(\hat x - x^*)\\
    &\Ra \delta\le -(\hat x+x^*)^T(\hat x - x^*) / t \le 4Rr/t.
\end{align}
Theorem \ref{thm:bound_opt} suggests that, as we push $t$ to $+\infty$ and solve a single optimization problem (\ref{eq:net}), we can find the optimizer of (\ref{eq:opt_net}). 

The remaining headache is the magnitude of $t$, say if it's too big, it may still be hard to implement (in fact, it's not a problem for today's machine accuracy). But we claim that, $t$ does not have to be crazy big. 

We review that, in \cite{chen2020learning}, the authors finds \emph{active constraints at optimum}, which is defined as the set of constraints that takes equality sign. They assume that the number of such active constraints at optimum is $n$, same as the dimension of the problem, and all active constraints are linearly independent. Then the neural network solves a combinatorial problem, which selects the active constraints and then the optimizer can be solved by linear equations. 

So, mirroring to our algorithm, if we can choose $t$ big enough, that only $n$ of the constraints are close enough to activation (if the upper and lower bound of the box are not too close), then we can identify the active constraints and then solve by linear equations. The overall algorithm is Algorithm \ref{algo2}.
\begin{algorithm}[htbp!]
\caption{Algorithm for solving (\ref{eq:opt_net})}\label{algo2}
\begin{algorithmic}
\REQUIRE{Parameters $c,u,v,A,b$, initial point $x_0$, accuracy $\delta$.}\\
\STATE{Compute radius $r$ and $R$. They can be upper bounded by the corresponding radius of box constraints (without linear equality constraints).}\vspace{0.3em}

\STATE{Set $t \la 4Rr/\delta$, $y\la -tc$.}\vspace{0.3em}

\STATE{Solve (\ref{eq:mu_prob}) by zeroth-order/gradient based algorithms. This is an unconstrained problem. Let the optimizer be $\mu^*$.}\vspace{0.3em}

\STATE{$(x^*,\lb_1^*,\lb_2^*) \la \phi(\mu^*)$.}
\IF{Return optimizer}
\RETURN{$x^*$.}
\ENDIF
\IF{Obtain active constraints and then solve linear equations}
\STATE{Select the active constraints by plugging $x^*$ in and pick the ones that are closest to equality.}\vspace{0.3em}

\STATE{Solve for optimizer $x_{\mathrm{les}}^*$ by applying the linear equation solver to active constraints.}
\RETURN{$x_{\mathrm{les}}^*$.}
\ENDIF
\end{algorithmic}
\end{algorithm}
\section{Acknowledgement}
The author would thank Yize Chen for the inspiring discussion of the model and technique in \cite{chen2020learning}. Finger crossed for his rebuttal.
\clearpage
\bibliographystyle{icml2019}
\bibliography{Bibfile}

\end{document}